\documentclass[12pt,a4paper,twoside]{article}

\usepackage[utf8]{inputenc}
\usepackage{graphicx}
\usepackage{tabularx}
\usepackage{textcomp}
\usepackage{kzartic}
\usepackage{caption}
\usepackage{hyperref}
\usepackage{commath}

\usepackage{natbib}
\usepackage{amsfonts}
\usepackage{amssymb}
\usepackage{amsthm}
\usepackage[tight,footnotesize]{subfigure}

\usepackage{tikz}
\usepackage{verbatim}
\usepackage{algorithm}% http://ctan.org/pkg/algorithms
\usepackage{lipsum}
\usepackage{algorithmic}

\newtheorem{definition}{Definition}
\newtheorem{theorem}{Theorem}

\newtheorem{example}{Example}
\newtheorem{remark}{Remark}

\DeclareCaptionLabelFormat{flushrightlabel}%
{\rightline{Table{~}#2}}
\DeclareCaptionLabelFormat{ryskk}%
{{Fig.~#2.}}
\usepackage{amsfonts}
\usepackage{times}
\renewcommand{\normalsize}{\usefont{T1}{ptm}{m}{n}\fontsize{13pt}{15pt}\selectfont}
\renewcommand{\bf}{\usefont{T1}{ptm}{b}{n}\fontsize{13pt}{15pt}\selectfont}

\renewcommand{\large}{\usefont{T1}{ptm}{m}{n}\fontsize{14pt}{16.8pt}\selectfont}

\usepackage{fancyhdr}
\fancypagestyle{plain}{%
	\fancyhead{}
	
}

\begin{document}
\global\count0=1
\setcounter{page}{1}
\makeatletter
\renewcommand{\@seccntformat }[1]{\csname the#1\endcsname\ }
\makeatother
\renewcommand{\textfraction}{0.0}
\renewcommand{\topfraction}{0.95}
\renewcommand{\bottomfraction}{0.95}
\renewcommand{\floatpagefraction}{0.9}
\thispagestyle{empty}
\vspace*{-1.75cm}
\noindent\underline{\hbox to \hsize{AUTOMATYZACJA PROCESÓW DYSKRETNYCH\hfill 2014}}\par
\vspace{2.2cm}
\sloppy\raggedbottom
\noindent  Konrad Andrzej MARKOWSKI\par
\noindent Warsaw University of Technology, Electrical Department\\
Institute of Control and Industrial Electronics\\
Koszykowa 75, 00-662 Warsaw\\
e-mail: Konrad.Markowski@ee.pw.edu.pl
\par
\vspace{1.45cm}
{\noindent\bf GRAFOWA METODA WYZNACZANIA MINIMALNEJ MACIERZY STANU UKŁADU UŁAMKOWEGO RZĘDU}\par

\vspace*{0.45cm}
\mbox{}\hfill\parbox[t]{14.75cm}{\renewcommand{\baselinestretch}{1}\footnotesize\normalsize%
\indent{\bf Streszczenie.}  
W pracy przedstawiono nową metodę wyznaczania macierzy stanu. Zaproponowana metoda bazuje na teorii grafów skierowanych. Przedstawiono procedurę wyznaczania elementów macierzy stanu a następnie zilustrowano ją przykładem numerycznym. }

\vspace*{0.6cm}
\noindent{\bf DIGRAPH METHOD FOR DETERMINATION OF MINIMAL STATE MATRIX OF THE FRACTIONAL SYSTEM}\par

\vspace*{0.45cm}
\mbox{}\hfill\parbox[t]{14.75cm}{\renewcommand{\baselinestretch}{1}\footnotesize\normalsize%
\indent{\bf Summary.}
A new method of determination entries of the state matrices has been presented. The presented method is based on one-dimensional digraph theory. A procedure for computation of the state matrices has also been proposed. The procedure has been illustrated with a numerical example.}

\renewcommand{\baselinestretch}{1}\footnotesize\normalsize
\vspace*{0.1cm}

\renewcommand{\recenzent}{Dr hab.~in.~Cezary Recenzent, prof.~Pol.~PL.}

\section{Introduction}
In the recent years many researchers were interested in positive linear systems \cite{Benvenuti2004}, \cite{Farina2000}, \cite{Kaczorek1985}, \cite{Kaczorek2003}, \cite{Luenberger1979}. In positive systems inputs, state variables and outputs take only non-negative values \cite{Berman1989}. Positive linear systems are defined on cones and not on linear spaces. Therefore, the theory of positive systems is more complicated than standard systems. The realisation problem is a very difficult task. In many research studies we can find canonical form of the system, i.e. constant matrix form, which satisfies the system described by the transfer function. With the use of this form we are able to write only one realisation of the system. In general we have a lot of solutions. This means that we can find many sets of matrices which fit into system transfer function. The state of the art in positive systems theory is given in the monographs \cite{Farina2000}, \cite{Kaczorek2003}, \cite{Luenberger1979}. A new method of determination entries of the state matrices will be proposed. A procedure for computation of the state matrices will be given. The procedure will be illustrated with a numerical example. 

This work has been organized as follows: Chapter 2 presents some notations and basic definitions of positive fractional systems and digraphs theory. In Chapter 3, we construct and discuss method for determination of the set of polynomial realisations which are based on digraphs theory. In Chapter 4, we illustrate method from Chapter 3 by numerical example. Finally, we give some concluding remarks, present open problems and bibliography positions.
\section{Preliminaries and Problem Formulation}

\subsection{Fractional System}
Consider the fractional discrete-time linear system, described by the equations \cite{Kaczorek2011}:
\begin{eqnarray}
\Delta ^{\alpha}x_{k+1}&=&\mathbf{A}x_{k}+\mathbf{B}u_{k}\nonumber\\
y_{k}&=&\mathbf{C}x_{k}+\mathbf{D}u_{k},\quad k\in\mathbb{Z}_{+}
\end{eqnarray}
where $x_{k}\in\mathbb{R}^{n}$, $u_{k}\in\mathbb{R}^{m}$, $y_{k}\in\mathbb{R}^{p}$ are the state, input and output vectors respectively and $\mathbf{A}\in\mathbb{R}^{n \times n}$, $\mathbf{B}\in\mathbb{R}^{n \times m}$, $\mathbf{C}\in\mathbb{R}^{p \times n}$, $\mathbf{D}\in\mathbb{R}^{p \times m}$, $\alpha\in\mathbf{R}$.

\begin{definition}\label{def:fractional}
The discrete-time function
\begin{eqnarray}
\Delta^{\alpha}x_{k}=\sum_{j=1}^{k}(-1)^{j}
	\left(
		\begin{array}{c}
		\alpha\\
		j
		\end{array}
	\right)x_{k-j}
\end{eqnarray}
where $0 < \alpha ,1$, $\alpha\in\mathbb{R}$ and
\begin{eqnarray}
\left(
		\begin{array}{c}
		\alpha\\
		j
		\end{array}
	\right)
=\left\{
\begin{array}{ll}
1 & \textit{for }k=0\\
\dfrac{(\alpha (\alpha -1)\dots (\alpha -k-1)}{k!} & \textit{for } k=1,2,3,\dots
\end{array}
\right.
\end{eqnarray}
is called the fractional $\alpha$ order difference of the function $x_{k}$.
\end{definition}
Using Definition \ref{def:fractional} we may write equation in the form:
\begin{eqnarray}\label{eq:FractionalSystemNew}
x_{k+1}+\sum_{j=1}^{k+1}(-1)^{j} \left(\begin{array}{c}\alpha \\ j \end{array}\right)x_{k-j+1}&=&\mathbf{A}x_{k}+\mathbf{B}u_{k}\\
y_{k}&=&\mathbf{C}x_{k}+\mathbf{D}u_{k}, \qquad k\in\mathbb{Z}_{+}\nonumber
\end{eqnarray}
\begin{definition}
The system (\ref{eq:FractionalSystemNew}) is called internally positive fractional system if $x_{k}\in\mathbb{R}_{+}^{n}$ and $y_{k}\in\mathbb{R}_{+}^{p}$, $k\in\mathbb{Z}_{+}$ for any initial conditions $x_{0}\in\mathbb{R}_{+}^{n}$ and all input sequence $u_{k}\in\mathbb{R}_{+}^{m}$, $k\in\mathbb{Z}_{+}$.
\end{definition}
\begin{definition}\label{def:PositiveMatrix}
The fractional system (\ref{eq:FractionalSystemNew}) is positive if and only if 
\begin{equation}\label{eq:PositiveMatrix}
\mathbf{A}+\alpha\mathbf{I}_{n}\in\mathbb{R}_{+}^{n \times n},\quad 0 < \alpha < 1,\quad 
\mathbf{B}\in\mathbb{R}_{+}^{n \times m}, \quad
\mathbf{C}\in\mathbb{R}_{+}^{p \times n}, \quad
\mathbf{D}\in\mathbb{R}_{+}^{p \times m}.
\end{equation}
\end{definition}
%
% ****************************************************************************************************
\subsection{Digraphs}
A directed graph  (or just digraph) $\mathfrak{D}$ consists of a non-empty finite set $\mathbb{V}(\mathfrak{D})$ of elements called vertices and a finite set $\mathbb{A}(\mathfrak{D})$ of ordered pairs of distinct vertices called arcs (\cite{Bang2009}). We call $\mathbb{V}(\mathfrak{D})$ the vertex set and $\mathbb{A}(\mathfrak{D})$ the arc set of $\mathfrak{D}$. We will often write $\mathfrak{D} = (\mathbb{V},\mathbb{A})$ which means that $\mathbb{V}$ and $\mathbb{A}$ are the vertex set and arc set of $\mathfrak{D}$, respectively. The order of $\mathfrak{D}$ is the number of vertices in $\mathfrak{D}$. The size of $\mathfrak{D}$ is the number of arc in $\mathfrak{D}$. For an arc $(v_{1},v_{2})$ the first vertex $v_{1}$ is its tail and the second vertex $v_{2}$ is its head.

There exists $\mathfrak{X}$-arc from vertex $v_{j}$ to vertex $v_{i}$ if and only if the $(i,j)$-th entry of the matrix $\mathbf{X}$ is nonzero. There exists $\mathfrak{Y}$ -arc from source $s_{l}$ to vertex $v_{j}$ if and only if the $l$-th entry of the matrix $\mathbf{Y}$ is nonzero.

\begin{example}
The system described by the following matrices
\begin{eqnarray}\label{eq:example_1}
	\left(\mathbf{X},\mathbf{Y}\right)=
	\left(
		\left[\begin{array}{ccc}
			1 & 0 & 1\\
			1 & 0 & 1\\
			1 & 1 & 0\\		
		\end{array}\right],					
		\left[\begin{array}{cc}
			1 & 0 \\
			0 & 0 \\
			0 & 1 \\			
		\end{array}\right],			
	\right)\nonumber
\end{eqnarray}
we can draw digraph $\mathfrak{D}^{(1)}$ consisting of vertices $v_{1}, v_{2}, v_{3}$ and source $s_{1}, s_{2}$. One-dimensional digraph corresponding to system (\ref{eq:example_1}) is presented in Figure \ref{fig:digraph_example_1}. 
\end{example}
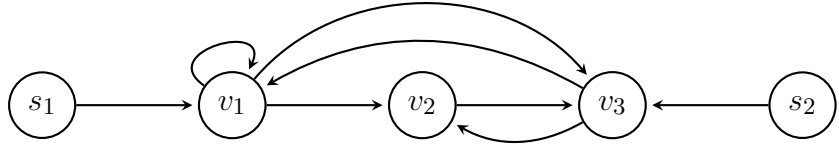
\begin{figure}[h]
\centering
\begin{tikzpicture}[scale=0.4, >=stealth,shorten >=2pt,node distance=2.5cm, thick,main node/.style={circle,draw}]
\node [main node] (v1) {$v_{1}$};
\node [main node] (v2) [right of=v1]{$v_{2}$};
\node [main node] (v3) [right of=v2]{$v_{3}$};
\node [main node] (s1) [left of=v1]{$s_{1}$};
\node [main node] (s2) [right of=v3]{$s_{2}$};
%\node [main node] (v2) [below right of=v1]{$v_{2}$};
%\node [main node] (v3) [above right of=v2]{$v_{3}$};

\path [->] 					(v1) edge node {} (v2)					%OK
							(v2) edge node {} (v3)					%OK
 							(v3) edge [bend right] node {} (v1)	%OK
 							(s1) edge node {} (v1)					%OK
 							(s2) edge node {} (v3)					%OK
				 			(v1) edge [out=145, in=65 , looseness=1.0, loop , distance=2cm] node {}  (v1) 		%OK
							(v1) edge [bend left=50] node {} (v3)	%OK
%							(v3) edge [bend right=40] node {} (v2)
%							(s1) edge node {} (v2)
%							(v1) edge [bend right] node {} (v2)
							(v3) edge [bend left] node  {} (v2);		%OK
%							(v3) edge [loop right] node {} (v3)
%							(s2) edge node {} (v2);

%\path [dashed, ->]	(v3) edge [bend right] node {} (v2);
\end{tikzpicture}
\caption{One-dimensional digraph }
\label{fig:digraph_example_1}
\end{figure}
We present below some basic notions from graph theory which are used in the farther considerations \cite{Bang2009}, \cite{Wallis2007}.

A walk in a digraphs $\mathfrak{D}^{2}$ is a finite sequence of arcs in which every two vertices $v_{i}$ and $v_{j}$ are adjacent or identical. For instance in Figure \ref{fig:digraph_example_1} there is the following walk: $(v_1 ,v_2)$, $(v_2, v_3)$, $(v_3,v_2)$, $(v_2 ,v_3 )$, $(v_3, v_1)$. A walk in which all of the arcs are distinct is called a path. For example in Figure \ref{fig:digraph_example_1} there is the following path: $(v_3, v_1)$, $(v_1,v_1)$, $(v_1,v_2)$, $(v_2,v_3)$, $(v_3,v_2)$. The path that goes through all vertices is called a finite path. For example in Figure \ref{fig:digraph_example_1} there is the following finite path: $(v_1,v_2)$, $(v_2,v_3)$. If the initial and the terminal vertices of the path are the same, then the path is called a cycle. For example in Figure \ref{fig:digraph_example_1} there are the following cycles: $(v_1,v_2)$, $(v_2,v_3)$, $(v_3,v_1)$.% , or $(v_2,v_3)$ , $(v_3,v_2)$.
%\newline

More information about use digraph theory in positive system is given in \cite{Fornasini1997}, \cite{Fornasini2005}, \cite{Markowski2013}, \cite{Markowski2014b}.
\subsection{Problem Formulation}
Formulation of the realisation problem for single-input single-output (SISO) positive fractional discrete time system is based on the following theorem.
\begin{theorem}
The transfer function of the fractional system (\ref{eq:FractionalSystemNew}) has the form:
\begin{eqnarray}\label{eq:TransferFunction}
T(z)=\mathbf{C}\left[ \mathbf{I}_{n}(z-c_{\alpha})-\mathbf{A} \right]^{-1}\mathbf{B}+\mathbf{D}
\end{eqnarray}
where:
\begin{eqnarray}
c_{\alpha}=c_{\alpha}(k-z)&=&\sum_{j=1}^{k+1}(-1)^{j-1}\left( \begin{array}{c} \alpha \\ j \end{array}\right)z^{1-j}\nonumber\\
N(z)&=&\mathbf{C}\textnormal{Adj}\left[ \mathbf{I}_{n}(z-c_{\alpha})-\mathbf{A} \right]+\mathbf{D}d(z) \nonumber\\
d(z)&=&\det\left[ \mathbf{I}_{n}(z-c_{\alpha})-\mathbf{A}\right]=
\nonumber
\end{eqnarray}
\end{theorem}
The proof of the Theorem 1 is given in \citep{Kaczorek2011}.

Matrices $\mathbf{A}$, $\mathbf{B}$, $\mathbf{C}$ and $\mathbf{D}$ satisfying (\ref{eq:PositiveMatrix}) are called a positive fractional realisation of a given transfer function described by the equation (\ref{eq:TransferFunction}). A realisation is called minimal if the dimension of the state matrix $\mathbf{A}$ is minimal among all realisations of $T(z)$.
\newline
\textbf{Our task is the following:} for a given transfer function (\ref{eq:TransferFunction}) determine entries of matrix $\mathbf{A}$ of the fractional system (\ref{eq:FractionalSystemNew}) using digraph $\mathfrak{D}^{(1)}$ theory. The dimension of the state matrices must be the minimal among possible ones and cannot appear additional condition on the coefficients of the characteristic polynomial. The problem of finding all possible realisations of a given transfer function is of such complexity that it cannot be solved in a reasonable time even by brute-force GPGPU method.
%
% **************************************************************************************
\section{Solution of The Problem}
%
% **************************************************************************************
%\subsection{Theorem}
%
By multiplying the numerator and denominator of (\ref{eq:TransferFunction}) by $(z-c_{\alpha})^{-n}$ we obtain:
\begin{eqnarray}\label{eq:TransferFunctionNew}
T(z^{-1})=\dfrac{b_{n}+b_{n-1}(z-c_{\alpha})^{-1}+\dots+b_{1}(z-c_{\alpha})^{1-n}+b_{0}(z-c_{\alpha})^{-n}}
{1-d_{n-1}(z-c_{\alpha})^{-1}-\dots-d_{1}(z-c_{\alpha})^{1-n}-d_{0}(z-c_{\alpha})^{-n}}
\end{eqnarray}
where 
\begin{eqnarray}
n(z^{-1})=b_{n}+b_{n-1}(z-c_{\alpha})^{-1}+\dots+b_{1}(z-c_{\alpha})^{1-n}+b_{0}(z-c_{\alpha})^{-n}
\end{eqnarray}
and
\begin{eqnarray}\label{eq:CharacteristicPolynomialNew}
d(z^{-1})=1-d_{n-1}(z-c_{\alpha})^{-1}-\dots-d_{1}(z-c_{\alpha})^{1-n}-d_{0}(z-c_{\alpha})^{-n}
\end{eqnarray}
is the characteristic polynomial.

Proposed method finds state matrix $\mathbf{A}$ using decomposing characteristic polynomial (\ref{eq:CharacteristicPolynomialNew}) into a set of simple monomials. 

In the first step we  decompose polynomial (\ref{eq:CharacteristicPolynomialNew}) into a set of the simple monomials in the following way:
\begin{eqnarray}\label{eq:CharacteristicPolynomialDecomposition}
d(z^{-1})=1-d_{n-1}(z-c_{\alpha})^{-1}-\dots-d_{1}(z-c_{\alpha})^{1-n}-d_{0}(z-c_{\alpha})^{-n}
\end{eqnarray}
For each simple monomials (\ref{eq:CharacteristicPolynomialDecomposition}) we create a digraph representation. Then we can determine all possible characteristic polynomial (\ref{eq:CharacteristicPolynomialNew}) realisations using all combinations of the digraph monomial representation. Finally, we combine received digraphs in one digraph which is corresponding to characteristic polynomial (\ref{eq:CharacteristicPolynomialNew}). 
 
\begin{theorem}\label{th:PropositionTheorem}
There exists positive state matrices $\mathbf{A}$ of the fractional discrete time linear system (\ref{eq:FractionalSystemNew}) corresponding to the characteristic polynomial (\ref{eq:CharacteristicPolynomialNew}) if
\begin{enumerate}
\item the coefficients of the characteristic polynomial
\begin{eqnarray}\label{eq:PropositionCondition1}
d_{i} \geqslant 0, \; for\; i=0,1,\dots n,\;\;d_{n}=1
\end{eqnarray}
\item the obtained digraph does not appear additional cycles and disjoint union. 
\end{enumerate}
\end{theorem}
\begin{proof}[Proof]
\textit{Condition 1.} The first condition came from Definition \ref{def:PositiveMatrix} and must be satisfied if we consider positive systems. If coefficients of the characteristic polynomial are negative then in state matrix negative elements appear. \newline
\textit{Condition 2.} Each monomial is represented by one cycle. If after combining all digraphs, which correspond to simple monomial, we obtain additional cycle this means that in polynomial additional simple monomial appears. If in digraph disjoint union appears this means that we do not have common parts for digraph corresponding to monomial and in polynomial additional simple monomial appears. In this situations we obtain a new polynomial which does not represent characteristic polynomial. 
\end{proof}
Using Theorem \ref{th:PropositionTheorem} we can construct $DetermineStateMatrix()$  algorithm. 
\begin{algorithm}[H]
\caption{$DetermineState Matrix()$}
\label{Alg:Main}
\begin{algorithmic}[1]
\algsetup{indent=2em}
\STATE $monomial = 1$; 
\STATE Determine number of $cycles$ in characteristic polynomial;
\FOR {$monomial = 1$ \TO {$cycles$}}
	\STATE Determine digraph $\mathfrak{D}^{(1)}$ for all monomial;
	\STATE $MonomialRealisation(monomial)$;
\ENDFOR	
\FOR {$monomial$ = 1 \TO {$cycles$}}
	\STATE Determine digraph as a combination of the digraph monomial representation
	\STATE $PolynomialRealisation(monomial)$;
	\IF {$PolynomialRealisation\;!=\;cycles$}
		\STATE Digraph contains additional cycles or digraph contains disjoint union
		\STATE \textbf{BREAK}
	\ELSIF{$PolynomialRealisation\;==\;cycles$}
		\STATE Digraph satisfies characteristic polynomial;
		\STATE Determine weights of the arcs in digraph; 
		\STATE Write state matrix $\mathbf{A}$;
		\RETURN $(PolynomialRealisation,\; \mathbf{A})$;
	\ENDIF
\ENDFOR
\end{algorithmic}
\end{algorithm}

%
% **************************************************************************************
\section{Numerical Example}
Let be given the transfer function 
\begin{eqnarray}\label{eq:Example}
T(z)=\dfrac{2(z-c_{\alpha})^{2}+5(z-c_{\alpha})+2}
{(z-c_{\alpha})^{2}-(z-c_{\alpha})-2}
\end{eqnarray}
Determine minimal state matrix which satisfies characteristic polynomial. 
By multiplying the numerator and denominator of (\ref{eq:Example}) by $(z-c_{\alpha})^{-2}$ we obtain transfer function in the form:
\begin{eqnarray}\label{eq:ExampleNew}
T(z^{-1})=\dfrac{2+5(z-c_{\alpha})^{-1}+2(z-c_{\alpha})^{-2}}{1-(z-c_{\alpha})^{-1}-2(z-c_{\alpha})^{-2}}
\end{eqnarray}
In the next step from transfer function (\ref{eq:ExampleNew}) we determine matrix $\mathbf{D}$ in the form:
\begin{eqnarray}
\mathbf{D}&=&\lim_{z\rightarrow \infty}T(z^{-1})=2\nonumber\\
T_{sp}&=&T(z^{-1})-\mathbf{D}=\dfrac{7(z-c_{\alpha})^{-1}+6(z-c_{\alpha})^{-2}}
{1-(z-c_{\alpha})^{-1}-2(z-c_{\alpha})^{-2}}
\end{eqnarray}
where
\begin{eqnarray}\label{eq:ExampleCharacteristicPolynomial}
d(z^{-1})=1-(z-c_{\alpha})^{-1}-2(z-c_{\alpha})^{-2}
\end{eqnarray}
is characteristic polynomial.

To determine all monomial realisation of the polynomial (\ref{eq:ExampleCharacteristicPolynomial}) in the first step we must write boundary conditions: 
\begin{itemize}
	\item number of vertices -- $VN=2$;
	\item number of colors in digraph -- $CN=1$;
	\item monomial: $M_{1}=(z-c_{\alpha})^{-1}$ and $M_{2}=2(z-c_{\alpha})^{-2}$;
	\item number of cycles -- $cycles=2$.
\end{itemize}

\begin{figure}[!h]
\centering
\begin{tikzpicture}[scale=0.25, >=stealth,shorten >=1pt,node distance=2.5cm, thick,main node/.style={circle,draw}]
\node [main node] (0) {$0$};
\node [main node] (1) [right of=0]{$1$};
\node [main node] (2) [right of=1]{$2$};
\path [->] 		(0) edge [in=100,out=80, distance=2.0cm] node [below]{\small $V[1,1]$} (1)
				(0) edge [in=100,out=85, distance=3.5cm] node [pos=0.72, below] {\small $V[1,2]$} (2);
\end{tikzpicture}
\caption{One-dimensional digraphs}
\label{fig:Example_1}
\end{figure}
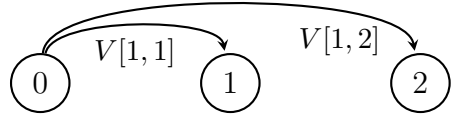

In the first step, of the algorithm we determine set of the possible connections between all vertices. Digraph $\mathfrak{D}^{(1)}$ presented in Figure \ref{fig:Example_1}.
\newline In the next step we determine all the possible realisations of the monomial $M_{1}$
\begin{figure}[!ht]
\centering
\subfigure[]{
\begin{tikzpicture}[scale=0.25, >=stealth,shorten >=1pt,node distance=2.5cm, thick,main node/.style={circle,draw}]
\node [main node] (0) {$0$};
\node [main node] (1) [right of=0]{$1$};
\node [main node] (2) [right of=1]{$2$};
\path [->] 					(0) edge node {} (1)
							(1) edge [out=145, in=45 , looseness=1.5, loop , distance=4cm] node [above]{\small $w(1,1)(z-c_{\alpha})^{-1}$}  (1);
\end{tikzpicture}
\label{fig:Example_2a}
}\quad
\subfigure[]{
\begin{tikzpicture}[scale=0.25, >=stealth,shorten >=1pt,node distance=2.5cm, thick,main node/.style={circle,draw}]
\node [main node] (0) {$0$};
\node [main node] (1) [right of=0]{$1$};
\node [main node] (2) [right of=1]{$2$};
\path [->] 					(0) edge node {} (1)
							(2) edge [out=145, in=45 , looseness=1.5, loop , distance=4cm] node [above] {\small $w(2,2)(z-c_{\alpha})^{-1}$} (2);
\end{tikzpicture}
\label{fig:Example_2b}
}
\caption{One-dimensional digraphs with all possible realisation of the monomial $M_{1}$}
\label{fig:Example_2}
\end{figure}
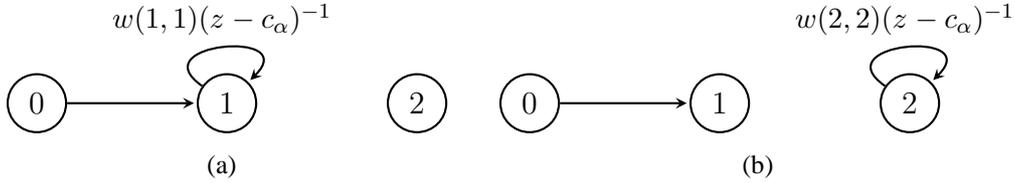
\newline In the same way we follow with monomial $M_{2}$. Digraph with all possible realisations of the monomial $M_{2}$ are presented in Figure 4.
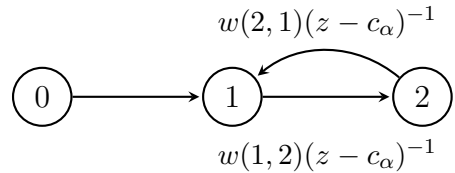
\begin{figure}[!ht]
\centering
\begin{tikzpicture}[scale=0.25, >=stealth,shorten >=1pt,node distance=2.5cm, thick,main node/.style={circle,draw}]
\node [main node] (0) {$0$};
\node [main node] (1) [right of=0]{$1$};
\node [main node] (2) [right of=1]{$2$};
\path [->] 					(0) edge node {} (1)
							(1) edge node [below, yshift=-0.4cm] {\small $w(1,2)(z-c_{\alpha})^{-1}$} (2)
							(2) edge [bend right=40] node [above] {\small $w(2,1)(z-c_{\alpha})^{-1}$} (1);
\end{tikzpicture}
\caption{One-dimensional digraph with all possible realisation of the monomial $M_{2}$}
\label{fig:Example_3}
\end{figure}
\newline In the next step we determine all combinations of the digraph monomial representation and write matrix $\mathbf{A}$:
\begin{enumerate}
\item Monomial $M_{1}$ from Figure \ref{fig:Example_2a} and monomial $M_{2}$ from Figure \ref{fig:Example_3}. Presented digraph (Figure \ref{fig:Example_4}) contains only two cycles and does not contain disjoint union. From digraph we can write set of the equation (Table \ref{tab:Table_1}). After solving them, we obtain weights coefficients and we can write matrix $\mathbf{A}$ in the following form:
\begin{eqnarray}
\mathbf{A}=\left[
\begin{array}{cc}
w(1,1) & w(2,1)\\
w(1,2) & 0
\end{array}
\right]
\end{eqnarray}
\begin{figure}[!ht]
\centering
\begin{tikzpicture}[scale=0.25, >=stealth,shorten >=1pt,node distance=2.5cm, thick,main node/.style={circle,draw}]
\node [main node] (0) {$0$};
\node [main node] (1) [right of=0]{$1$};
\node [main node] (2) [right of=1]{$2$};
\path [->] 		(0) edge node {} (1)
				(1) edge node [below, yshift=-0.4cm] {\small $w(1,2)(z-c_{\alpha})^{-1}$} (2)
				(2) edge [bend right=40] node [above, xshift=1.0cm] {\small $w(2,1)(z-c_{\alpha})^{-1}$} (1)
				(1) edge [out=145, in=65 , looseness=1.5, loop , distance=4cm] node [above, xshift=-1.0cm] {\small $w(1,1)(z-c_{\alpha})^{-1}$}  (1);
\end{tikzpicture}
\caption{Realisation of the characteristic polynomial (\ref{eq:ExampleCharacteristicPolynomial})}
\label{fig:Example_4}
\end{figure}
\begin{table}[!h]
\begin{center}
\caption{Set of the equations for digraph presented in Figure 5 }
\begin{tabular}{||c|c|c||}
\hline\hline
 Power & Digraph & Polynomial  \tabularnewline
 \hline\hline
  $(z-c_{\alpha})^{-1}$ & $w(1,1)$  & 1\tabularnewline
 \hline
  $(z-c_{\alpha})^{-2}$ & $w(1,2)w(2,1)$ & 2\tabularnewline
 \hline
 \hline
\end{tabular}
\label{tab:Table_1}
\end{center}
\end{table}
\item Monomial $M_{1}$ from Figure \ref{fig:Example_2b} and monomial $M_{2}$ from Figure \ref{fig:Example_3}. Presented digraph (Figure \ref{fig:Example_5}) contains only two cycles and does not contain disjoint union. From digraph we can write set of the equation (Table \ref{tab:Table_2}). After solving them, we obtain weights coefficients and we can write matrix $\mathbf{A}$ in the following form:
\begin{eqnarray}
\mathbf{A}=\left[
\begin{array}{cc}
0 & w(2,1)\\
w(1,2) & w(2,2)
\end{array}
\right]
\end{eqnarray}
\begin{figure}[!ht]
\centering
\begin{tikzpicture}[scale=0.25, >=stealth,shorten >=1pt,node distance=2.5cm, thick,main node/.style={circle,draw}]
\node [main node] (0) {$0$};
\node [main node] (1) [right of=0]{$1$};
\node [main node] (2) [right of=1]{$2$};
\path [->] 			(0) edge node {} (1)
					(1) edge node [below, yshift=-0.4cm] {\small $w(1,2)(z-c_{\alpha})^{-1}$} (2)
					(2) edge [bend right=40] node [above, xshift=-0.5cm] {\small $w(2,1)(z-c_{\alpha})^{-1}$} (1)
					(2) edge [out=125, in=45 , looseness=1.5, loop , distance=4cm] node [above, xshift=1.5cm] {\small $w(2,2)(z-c_{\alpha})^{-1}$}  (2);
\end{tikzpicture}
\caption{Realisation of the characteristic polynomial (\ref{eq:ExampleCharacteristicPolynomial})}
\label{fig:Example_5}
\end{figure}
\begin{table}[!h]
\begin{center}
\caption{Set of the equations for digraph presented in Figure 6 }
\begin{tabular}{||c|c|c||}
\hline\hline
 Power & Digraph & Polynomial  \tabularnewline
 \hline\hline
  $(z-c_{\alpha})^{-1}$ & $w(2,2)$  & 1\tabularnewline
 \hline
  $(z-c_{\alpha})^{-2}$ & $w(1,2)w(2,1)$ & 2\tabularnewline
 \hline
 \hline
\end{tabular}
\label{tab:Table_2}
\end{center}
\end{table}
\end{enumerate}
\begin{remark}
By solving the set of the equation presented in Table \ref{tab:Table_1} or Table \ref{tab:Table_2} we can determine a lot of state matrices $\mathbf{A}$ which satisfy its characteristic polynomial (\ref{eq:ExampleCharacteristicPolynomial}).
\end{remark}
%
% **************************************************************************************
\section{Concluding Remarks}
The paper includes a simple method based on digraph theory to determine minimal realisation of the characteristic polynomial of a positive one dimensional fractional system. By using this method, a fas algorithm for determining all possible realisations of the characteristic polynomial was constructed. The proposed algorithm is based on the digraphs theory. Currently, the method of determining a positive polynomial realisation using GPU units and digraphs methods is being implemented in the memory-efficient way. At the same time we are working on extension of the presented algorithm to solve reachability and realisation problems. Extending the proposed algorithm to dynamic systems of another class as well as searching for new areas of using multiprocessing calculations remains an open problem.
%
% **************************************************************************************

\bibliographystyle{plainnat}
\bibliography{kkapd}
%\begin{thebibliography}{9}
%\itemsep 3pt
%
%\bibitem{AAA}
%Nazwisko autora, inicjały imion: Tytuł. Wydawnictwo, nr tomu, miejsce wydania, rok, strony.
%\bibitem{BBB}
%Author A., Author B., Author C.: {Title}. Journal, 25(2), 1998, p.\ 115--145.
%\bibitem{CCC}
%Author D., Author E., Author F.G.: {Title}. Proc. 16th IEEE
%Symposium (2005), p.\ 440--445.
%
%\end{thebibliography}

\end{document}